\documentstyle[12pt,amssymb]{article}

\begin{document}

\baselineskip16pt

\newtheorem{definition}{Definition $\!\!$}[section]
\newtheorem{prop}[definition]{Proposition $\!\!$}
\newtheorem{lem}[definition]{Lemma $\!\!$}
\newtheorem{corollary}[definition]{Corollary $\!\!$}
\newtheorem{theorem}[definition]{Theorem $\!\!$}
\newtheorem{example}[definition]{\it Example $\!\!$}
\newtheorem{remark}[definition]{Remark $\!\!$}

\newcommand{\nc}[2]{\newcommand{#1}{#2}}
\newcommand{\rnc}[2]{\renewcommand{#1}{#2}}

\nc{\Section}{\setcounter{definition}{0}\section}
\renewcommand{\theequation}{\thesection.\arabic{equation}}
\newcounter{c}
\renewcommand{\[}{\setcounter{c}{1}$$}
\newcommand{\etyk}[1]{\vspace{-7.4mm}$$\begin{equation}\label{#1}
\addtocounter{c}{1}}
\renewcommand{\]}{\ifnum \value{c}=1 $$\else \end{equation}\fi}


\nc{\bpr}{\begin{prop}}
\nc{\bth}{\begin{theorem}}
\nc{\ble}{\begin{lem}}
\nc{\bco}{\begin{corollary}}
\nc{\bre}{\begin{remark}}
\nc{\bex}{\begin{example}}
\nc{\bde}{\begin{definition}}
\nc{\ede}{\end{definition}}
\nc{\epr}{\end{prop}}
\nc{\ethe}{\end{theorem}}
\nc{\ele}{\end{lem}}
\nc{\eco}{\end{corollary}}
\nc{\ere}{\hfill\mbox{$\Diamond$}\end{remark}}
\nc{\eex}{\end{example}}
\nc{\epf}{\hfill\mbox{$\Box$}}
\nc{\ot}{\otimes}
\nc{\bsb}{\begin{Sb}}
\nc{\esb}{\end{Sb}}
\nc{\ct}{\mbox{${\cal T}$}}
\nc{\ctb}{\mbox{${\cal T}\sb B$}}
\nc{\bcd}{\[\begin{CD}}
\nc{\ecd}{\end{CD}\]}
\nc{\ba}{\begin{array}}
\nc{\ea}{\end{array}}
\nc{\bea}{\begin{eqnarray}}
\nc{\eea}{\end{eqnarray}}
\nc{\be}{\begin{enumerate}}
\nc{\ee}{\end{enumerate}}
\nc{\beq}{\begin{equation}}
\nc{\eeq}{\end{equation}}
\nc{\bi}{\begin{itemize}}
\nc{\ei}{\end{itemize}}
\nc{\kr}{\mbox{Ker}}
\nc{\te}{\!\ot\!}
\nc{\pf}{\mbox{$P\!\sb F$}}
\nc{\pn}{\mbox{$P\!\sb\nu$}}
\nc{\bmlp}{\mbox{\boldmath$\left(\right.$}}
\nc{\bmrp}{\mbox{\boldmath$\left.\right)$}}
\rnc{\phi}{\mbox{$\varphi$}}
\nc{\LAblp}{\mbox{\LARGE\boldmath$($}}
\nc{\LAbrp}{\mbox{\LARGE\boldmath$)$}}
\nc{\Lblp}{\mbox{\Large\boldmath$($}}
\nc{\Lbrp}{\mbox{\Large\boldmath$)$}}
\nc{\lblp}{\mbox{\large\boldmath$($}}
\nc{\lbrp}{\mbox{\large\boldmath$)$}}
\nc{\blp}{\mbox{\boldmath$($}}
\nc{\brp}{\mbox{\boldmath$)$}}
\nc{\LAlp}{\mbox{\LARGE $($}}
\nc{\LArp}{\mbox{\LARGE $)$}}
\nc{\Llp}{\mbox{\Large $($}}
\nc{\Lrp}{\mbox{\Large $)$}}
\nc{\llp}{\mbox{\large $($}}
\nc{\lrp}{\mbox{\large $)$}}
\nc{\lbc}{\mbox{\Large\boldmath$,$}}
\nc{\lc}{\mbox{\Large$,$}}
\nc{\Lall}{\mbox{\Large$\forall\;$}}
\nc{\bc}{\mbox{\boldmath$,$}}
\rnc{\epsilon}{\varepsilon}
\rnc{\ker}{\mbox{\em Ker}}
\nc{\ra}{\rightarrow}
\nc{\ci}{\circ}
\nc{\cc}{\!\ci\!}
\nc{\T}{\mbox{\sf T}}
\nc{\can}{\mbox{\em\sf T}\!\sb R}
\nc{\cnl}{$\mbox{\sf T}\!\sb R$}
\nc{\lra}{\longrightarrow}
\nc{\M}{\mbox{Map}}
\rnc{\to}{\mapsto}
\nc{\imp}{\Rightarrow}
\rnc{\iff}{\Leftrightarrow}
\nc{\ob}{\mbox{$\Omega\sp{1}\! (\! B)$}}
\nc{\op}{\mbox{$\Omega\sp{1}\! (\! P)$}}
\nc{\oa}{\mbox{$\Omega\sp{1}\! (\! A)$}}
\nc{\inc}{\mbox{$\,\subseteq\;$}}
\nc{\de}{\mbox{$\Delta$}}
\nc{\spp}{\mbox{${\cal S}{\cal P}(P)$}}
\nc{\dr}{\mbox{$\Delta_{R}$}}
\nc{\dsr}{\mbox{$\Delta_{\cal R}$}}
\nc{\m}{\mbox{m}}
\nc{\0}{\sb{(0)}}
\nc{\1}{\sb{(1)}}
\nc{\2}{\sb{(2)}}
\nc{\3}{\sb{(3)}}
\nc{\4}{\sb{(4)}}
\nc{\5}{\sb{(5)}}
\nc{\6}{\sb{(6)}}
\nc{\7}{\sb{(7)}}
\nc{\hsp}{\hspace*}
\nc{\nin}{\mbox{$n\in\{ 0\}\!\cup\!{\Bbb N}$}}
\nc{\al}{\mbox{$\alpha$}}
\nc{\bet}{\mbox{$\beta$}}
\nc{\ha}{\mbox{$\alpha$}}
\nc{\hb}{\mbox{$\beta$}}
\nc{\hg}{\mbox{$\gamma$}}
\nc{\hd}{\mbox{$\delta$}}
\nc{\he}{\mbox{$\varepsilon$}}
\nc{\hz}{\mbox{$\zeta$}}
\nc{\hs}{\mbox{$\sigma$}}
\nc{\hk}{\mbox{$\kappa$}}
\nc{\hm}{\mbox{$\mu$}}
\nc{\hn}{\mbox{$\nu$}}
\nc{\la}{\mbox{$\lambda$}}
\nc{\hl}{\mbox{$\lambda$}}
\nc{\hG}{\mbox{$\Gamma$}}
\nc{\hD}{\mbox{$\Delta$}}
\nc{\th}{\mbox{$\theta$}}
\nc{\Th}{\mbox{$\Theta$}}
\nc{\ho}{\mbox{$\omega$}}
\nc{\hO}{\mbox{$\Omega$}}
\nc{\hp}{\mbox{$\pi$}}
\nc{\hP}{\mbox{$\Pi$}}
\nc{\bpf}{{\it Proof.~~}}
\nc{\slq}{\mbox{$A(SL\sb q(2))$}}
\nc{\fr}{\mbox{$Fr\llp A(SL(2,\IC))\lrp$}}
\nc{\slc}{\mbox{$A(SL(2,\IC))$}}
\nc{\af}{\mbox{$A(F)$}}
\rnc{\widetilde}{\tilde}
\nc{\suq}{\mbox{$A(SU_q(2))$}}
\nc{\asq}{\mbox{$A(S_q^2)$}}
\nc{\tasq}{\mbox{$\widetilde{A}(S_q^2)$}}

\def\esl{{\mbox{$E\sb{\frak s\frak l (2,{\Bbb C})}$}}}
\def\esu{{\mbox{$E\sb{\frak s\frak u(2)}$}}}
\def\zf{{\mbox{${\Bbb Z}\sb 4$}}}
\def\zt{{\mbox{$2{\Bbb Z}\sb 2$}}}
\def\ox{{\mbox{$\Omega\sp 1\sb{\frak M}X$}}}
\def\oxh{{\mbox{$\Omega\sp 1\sb{\frak M-hor}X$}}}
\def\oxs{{\mbox{$\Omega\sp 1\sb{\frak M-shor}X$}}}
\def\Fr{\mbox{Fr}}
\def\gal{-Galois extension}
\def\hge{Hopf-Galois extension}
\def\ta{\tilde a}
\def\tb{\tilde b}
\def\tc{\tilde c}
\def\td{\tilde d}
\def\st{\stackrel}


\newcommand{\Sp}{{\rm Sp}\,}
\newcommand{\Mor}{\mbox{$\rm Mor$}}
\newcommand{\skrA}{{\cal A}}
\newcommand{\Phase}{\mbox{$\rm Phase\,$}}
\newcommand{\id}{{\rm id}}
\newcommand{\diag}{{\rm diag}}
\newcommand{\inv}{{\rm inv}}
\newcommand{\ad}{{\rm ad}}
\newcommand{\poi}{{\rm pt}}
\newcommand{\Dim}{{\rm dim}\,}
\newcommand{\Ker}{{\rm ker}\,}
\newcommand{\Mat}{{\rm Mat}\,}
\newcommand{\Rep}{{\rm Rep}\,}
\newcommand{\Fun}{{\rm Fun}\,}
\newcommand{\Tr}{{\rm Tr}\,}
\newcommand{\supp}{\mbox{$\rm supp$}}
\newcommand{\half}{\frac{1}{2}}
\newcommand{\skrF}{{A}}
\newcommand{\skrD}{{\cal D}}
\newcommand{\skrC}{{\cal C}}
\newcommand{\ttimes}{\mbox{$\hspace{.5mm}\bigcirc\hspace{-4.9mm}
\perp\hspace{1mm}$}}
\newcommand{\Ttimes}{\mbox{$\hspace{.5mm}\bigcirc\hspace{-3.7mm}
\raisebox{-.7mm}{$\top$}\hspace{1mm}$}}
\newcommand{\bbr}{{\bf R}}
\newcommand{\bbz}{{\bf Z}}
\newcommand{\Ci}{C_{\infty}}
\newcommand{\Cb}{C_{b}}
\newcommand{\fa}{\forall}
\newcommand{\rrr}{right regular representation}
\newcommand{\wrt}{with respect to}
\newcommand{\qg}{quantum group}
\newcommand{\qgs}{quantum groups}
\newcommand{\cs}{classical space}
\newcommand{\qs}{quantum space}
\newcommand{\po}{Pontryagin}
\newcommand{\ch}{character}
\newcommand{\chs}{characters}

\def\inbar{\,\vrule height1.5ex width.4pt depth0pt}
\def\IC{{\Bbb C}}
\def\IZ{{\Bbb Z}}
\def\IN{{\Bbb N}}
\def\otc{\otimes_{\IC}}
\def\ra{\rightarrow}
\def\ota{\otimes_ A}
\def\otza{\otimes_{ Z(A)}}
\def\otc{\otimes_{\IC}}
\def\h{\rho}
\def\x{\zeta}
\def\th{\theta}
\def\s{\sigma}
\def\t{\tau}
\nc{\tens}{\otimes}
\nc{\aqpb}{algebraic quantum principal bundle}
\nc{\qpb}{quantum principal bundle}
\nc{\hop}{Hopf algebra}
\def\eps{{\epsilon}}
\def\extd{{\rm d}}
\def\<{\langle}
\def\>{\rangle}
\def\isom{{\cong}}
\def\Hom{{\rm Hom}}
\def\Ad{{\rm Ad}}
\def\ev{{\rm ev}}
\def\coev{{\rm coev}}
\def\id{{\rm id}}
\def\o{{}_{(1)}}\def\t{{}_{(2)}}\def\th{{}_{(3)}}
\def\bo{{}^{\bar{(1)}}}\def\bt{{}^{\bar{(2)}}}
\def\note#1{}
\def\qqquad{\qquad\quad}
\def\nquad{{\!\!\!\!\!\!}}
\def\nqquad{\nquad\nquad}\def\nqqquad{\nqquad\nquad}
\def\equad{\nquad}
\def\eqn#1#2{\begin{equation}#2\label{#1}\end{equation}}
\def\align#1{\begin{eqnarray*}#1\end{eqnarray*}}
\def\alignn#1#2{\begin{eqnarray}\label{#1}#2\end{eqnarray}}
\def\cmath#1{\[\begin{array}{c} #1 \end{array}\]}
\def\ceqn#1#2{\begin{equation}\label{#1}\begin{array}{c}#2\end{array}
\end{equation}}
\renewcommand\arraystretch{1.7}

\title{
\vspace*{-15mm}
{\LARGE\bf PROJECTIVE MODULE DESCRIPTION OF THE Q-MONOPOLE
\thanks{DAMTP-97-115 
}
}\\
}
\author{
{\normalsize\sc Piotr M.~Hajac}\thanks{
http://www.damtp.cam.ac.uk/user/pmh33  (E-mail: pmh33@amtp.cam.ac.uk)
On leave from:
Department of Mathematical Methods in Physics, 
Warsaw University, ul.~Ho\.{z}a 74, Warsaw, \mbox{00--682~Poland}.
http://info.fuw.edu.pl/KMMF/ludzie\underline{~~}ang.html 
(E-mail: pmh@fuw.edu.pl)} 
\normalsize ~and~ 
{\normalsize\sc Shahn Majid}\thanks{
http://www.damtp.cam.ac.uk/user/majid  (E-mail: majid@amtp.cam.ac.uk)}\\
\normalsize Department of Applied Mathematics and Theoretical Physics,
\vspace*{-5mm}\\
\normalsize University of Cambridge, Silver St., Cambridge CB3 9EW, England.\\
\vspace*{-15mm}
}
\date{}
\maketitle

\vspace*{-5mm}
\begin{abstract}
The Dirac $q$-monopole connection is used to compute projector matrices
of quantum Hopf line bundles for arbitrary winding number.
The Chern-Connes pairing of cyclic cohomology and $K$-theory is computed for 
 the winding number $-1$. The non-triviality of this pairing is used
to conclude that the quantum principal Hopf fibration is non-cleft.
 Among general results, we provide a left-right symmetric 
characterization of the canonical strong connections on 
quantum principal homogeneous spaces with an injective antipode.
We also provide for arbitrary strong connections on algebraic quantum principal 
bundles (Hopf-Galois extensions) their associated covariant derivatives 
on projective modules. 
\end{abstract}

\section*{Introduction}

The goal of this paper is to provide a better understanding of the relationship
between the quantum-group and $K$-theory approach to the noncommutative-geometry
gauge theory. The latter approach is based on the classical Serre-Swan theorem 
that allows one to think of vector bundles as projective modules. The former
comes from the concept of a Hopf-Galois extension which describes a quantum
principal bundle the same way Hopf algebras describe \qg s. Here
a Hopf algebra $H$ plays the role of the algebra of functions on the
structure group, and the total space of 
a bundle is replaced by an  $H$-comodule algebra~$P$. 
We rely on the Hopf-Galois
theory to derive our 
noncommutative-geometric constructions. On the other hand, it is the machinery
of noncommutative geometry that allows
us to obtain a Galois-theoretic result: We employ the Chern-Connes pairing 
to prove the non-cleftness of the \hge\
of the algebraic quantum principal Hopf fibration.

We begin in Section~1 with some preliminaries about \hge s, connections and 
connection 1-forms on \aqpb s, and connections on projective modules.
In Section~2 we extend the existing theory with some 
general results about strong connections, their covariant derivatives on
projective modules, and bicovariant splittings of canonical \hop\ surjections. 
We also discuss how to obtain projector matrices from splittings of 
the multiplication map. In
Section~3 we first define (the space of sections of) a quantum Hopf line bundle
as a bimodule associated to the quantum principal Hopf fibration 
via a one-dimensional corepresentation of the Hopf
algebra $k[z,z^{-1}]$. Then we use a canonical strong connection on 
the quantum principal Hopf fibration 
(Dirac $q$-monopole) to compute, for any one-dimensional
corepresentation, left and right projector matrices of the thus defined quantum
Hopf line bundles. This computation is the main part of our paper and
provides the projective-module characterization of the $q$-monopole.
Further results relating to the Chern-Connes pairing are in Section~4.
We end with Appendix where we show that the only invertible elements of the
coordinate ring of $SL_q(2)$ are non-zero numbers, and use it as an alternative way
to conclude the non-cleftness of the quantum Hopf fibration. 

To focus attention and 
take advantage of the cyclic cohomology results in~\cite{mnw91}, 
we work over a ground
field $k$ of characteristic zero, and assume that $q$ is a non-zero element in $k$
that is not a root of~1. 
We use the Sweedler notation $\Delta h=h\o\tens h\t$
(summation understood) and its derivatives. The antipode of the
Hopf algebra is a linear map $S:H\ra H$, and the counit is an algebra
map $\eps:H\ra k$ obeying certain properties. The convolution product of two
linear maps from a coalgebra to an algebra is denoted in the following
way: $(f*g)(c):=f(c\1)g(c\2)$. We use interchangeably the words ``colinear"
and ``covariant" with respect to linear maps that preserve the comodule structure.
For an introduction to noncommutative geometry, 
quantum groups, Hopf-Galois extensions and quantum-group gauge theory 
we refer to \cite{c-a94,l-g97}, 
\cite{m-s95}, \cite{s-hj94} and \cite{bm93,bm97} respectively.

\section{Preliminaries}

We begin by recalling basic definitions and known results.
\bde\label{pmdef}
Let $\cal E$ be a left \mbox{$B$-module}, and
$(\Omega(B),d)$ a differential algebra on~$B$. A linear map 
$\nabla : \Omega\sp*(B)\ot\sb B{\cal E}\ra\Omega\sp{*+1}(B)\ot\sb B{\cal E}$ 
is called a connection (covariant derivative) on $\cal E$ iff
$
{\forall}\;\xi\in{\cal E},\, \hl\in\Omega(\! B)\, :\; 
\nabla(\hl\ot\sb B\xi)=\hl(\nabla\xi) +d\hl\ot\sb B\xi\, . 
$
\ede
In the case of the universal differential algebra
the existence of a connection is {\em equivalent} to the projectivity 
of~$\cal E$~\cite[Corollary~8.2]{cq95}, \cite[Proposition~8.2.3]{l-g97}. 
If $\cal E$ is projective then
a connection exists for any differential algebra because it can be obtained
from the universal differential algebra and the canonical surjection onto
a given differential algebra~\cite[p.555]{c-a94}.
\bde\label{hgdef} 
Let $H$ be a Hopf algebra, $P$ be a right $H$-comodule algebra with multiplication
$m_P$ and coaction \dr, and 
$B:=P\sp{co H}:=\{p\in P\,|\; \dr\ p=p\ot 1\}$ the subalgebra of coinvariants. 
We say that P is a (right) {\em $H$-Galois extension} of $B$ iff the canonical 
left $P$-module right $H$-comodule map 
\[
\chi:=(m\sb P\ot id)\circ (id\ot\sb B \dr )\, :\; P\ot\sb B P\lra P\ot H
\]
is bijective. We say that $P$ is a faithfully flat $H$\gal\ of $B$
iff $P$ is faithfully flat as a right and left $B$-module.
\ede
For a comprehensive review of the concept of faithful flatness see~\cite{b-n72}.
\bde\label{cleft} 
An $H$\gal\ is called {\em cleft} iff there exists a unital convolution
invertible linear map $\Phi : H\ra P$ satisfying 
$\dr\circ\Phi=(\Phi\ot id)\circ\hD$. We call $\Phi$ a cleaving map of~$P$.
\ede
Note that, in general, $\Phi$ is {\em not} uniquely determined by its defining 
conditions. Observe also that the unitality assumption for the cleaving map
is unnecessary in the sense that any right colinear convolution invertible mapping
can be normalised to be unital. Indeed, let $\tilde{\Phi}$ be such a mapping,
and $\tilde{\Phi}(1):=b$. By the colinearity, we have that $b\in B$, and
 the convolution invertibility entails that $b$ is invertible. Also, 
$b^{-1}\ot 1=b^{-1}\dr(bb^{-1})=b^{-1}b\dr(b^{-1})=\dr(b^{-1})$. 
It is straightforward to check that $\Phi:=b^{-1}\tilde{\Phi}$ is right colinear,
 convolution invertible and unital. Let us also remark that a cleaving map
is necessarily injective:
\[
\llp m_P\ci(m_P\ot id)\ci(id\ot\Phi^{-1}\ot id)\ci(id\ot\hD)\ci\dr\ci\Phi\lrp(h)
=\Phi(h\1)\Phi^{-1}(h\2)h\3=h,~~~
\forall\, h\in H.
\]

To fix convention, let us recall that the universal differential calculus
 (grade one of the universal differential algebra) can be defined as the kernel of
 the multiplication map $\hO^1B:=\mbox{Ker}(B\ot B\st{m}{\ra}B)$ with the 
differential $db:=1\ot b-b\ot 1$ (e.g., see~\cite[Section~7.1]{l-g97}). 
(We abuse the notation and use the same letter
$d$ to signify both the universal and general differential.) 
The following are the universal-differential-calculus 
 versions of more general definitions in~\cite{bm93,h-pm96}:
\bde[\cite{bm93}]\label{condef}
Let $B\inc P$ be an $H$\gal. Denote by $\hO^1P$ the universal differential
calculus on~$P$. A left $P$-module projection $\Pi$ on $\hO^1P$ is 
called a {\em connection} on a \qpb~iff
\be
\item $\ker\,\Pi =  P(\Omega\sp{1}B)P$ (horizontal forms),
\item $\dsr\ci\Pi = (\Pi\ot id)\ci\dsr$ (right covariance).
\ee
\ede 
Here \dsr\ is the right coaction on differential forms given by the formula
$\dsr (ada'):=a\0 da'\0\ot a\1 a'\1$, where $\dr a:=a\0\ot a\1$ 
(summation understood).
Coaction on higher order forms is defined in the same manner.
\bde[\cite{bm93}] \label{confordef}
Let $P$, $H$, $B$ and $\hO^1P$ be as above.
 A $k$-homomorphism $\omega : H\ra\hO^1P$ such that $\ho(1)=0$
is called a
{\em connection form} iff it satisfies the
following properties:
\be
\item $ (m\sb P\ot id)\ci(id\ot\dr )\ci\omega = 1\ot(id - \eps)$ 
(fundamental vector field condition),
\item $\dsr\ci\omega = (\omega\ot id)\ci
ad_{R}$, $ad_{R}(h):=h\2\ot S(h\1)h\3$ (right adjoint covariance).
\ee
\ede
For every \hge\ there is a one-to-one
correspondence between connections and
connection forms (see~\cite[Proposition~2.1]{m-s97}). 
In particular, the connection
$\Pi\sp{\omega}$ associated to
a connection form $\omega$ is given by the
formula: 
\beq\label{omefor}
\Pi\sp{\omega}(dp)=p\0\ho(p\1)\, .
\eeq
$\Pi\sp{\omega}$ is a left $P$-module homomorphism, so that it suffices to know 
its values on exact forms.
\bde[\cite{h-pm96}]\label{scdef}
Let $\Pi$ be a connection in the sense of Definition~\ref{condef}. It
 is called {\em strong} iff 
$(id - \Pi)(dP)\inc(\Omega\sp{1}B)P$. We say that a connection form is strong
iff its associated connection is strong.
\ede

A natural next step is to consider
associated quantum vector bundles. More precisely, what we need here is a
replacement of the module of sections of an associated vector bundle. In the
classical case such sections can be equivalently described  as ``functions
of type $\varrho$" from the total space of a principal bundle to a vector 
space. We follow this construction in the quantum case by considering $B$-bimodules
of colinear maps $\mbox{Hom}_\rho(V,P)$ associated with 
an $H$\gal\ $B\inc P$ via a corepresentation $\rho:V\ra V\ot H$ (see~\cite{d-m96}). 
For our later purpose, we need the following reformulation of
\cite[Proposition~A.7]{bm93}:
\ble\label{isol}
Let $B\inc P$ be a cleft $H$\gal\ and 
$\rho:V\ra V\ot H$ a right corepresentation of $H$ on $V$.
Then the space of colinear maps
$\mbox{\em Hom}_\rho(V,P)$ is isomorphic as a left 
$B$-module  to the free module $\mbox{\em Hom}(V,B)$.
\ele

\section{Strong connections on associated projective modules}

First we study a general setting for translating strong
connections on algebraic quantum principal bundles to connections on
projective modules. The associated bimodule of colinear maps is 
finitely generated 
projective as a left module over the subalgebra of coinvariants under
rather unrestrictive assumptions. However, we do not assume the
projectivity of this module in the following two propositions, as it is
needed only later to ensure the existence of a connection. Also,
although we work only with the universal differential algebra in the sequel,
we do not assume here that the differential algebra is universal. It
suffices that it is right-covariant, i.e., the right coaction is
well-defined on differential forms, and right-covariant and right-flat
in the second proposition. On the other hand, we do not aim here at
the utmost generality but try to keep our noncommutative-geometric motivation
evident.

\bpr\label{forms}
Let $H$ be a Hopf algebra with a bijective antipode,
$P$ a faithfully flat $H$-Galois extension of $B$, and 
$V\st{\rho}{\ra}V\ot H$  (dim$V<\infty$) a coaction.
Denote by $\mbox{\em Hom}_\rho(V,P)$ the $B$-bimodule of colinear 
homomorphisms from $V$ to $P$, 
and choose a right-covariant differential algebra $\hO(P)$. 
Then the following map 
\[
\check{\ell}:\hO(B)\ot_B\mbox{\em Hom}_\rho(V,P)\lra
\mbox{\em Hom}_\rho(V,\hO(B)P),~~~
\llp\check{\ell}(\hl\ot_B\phi)\lrp(v)=\hl\phi(v),
\]
is an isomorphism of graded left $\hO(B)$-modules.
\epr\bpf
It suffices to show that $\check{\ell}$ has an inverse. By choosing a linear
basis $\{\hl_\mu\}$ of $\hO(B)$, for any $\phi\in\mbox{Hom}_\rho(V,\hO(B)P)\,$
 we can write $\phi(v)=\sum_\mu\hl_\mu\phi^\mu(v)$. The point is now to show
that we can always choose each $\phi^\mu$ to be an element of
$\mbox{Hom}_\rho(V,P)$. It can be done by assuming flatness of $\hO(B)$
(see Proposition~\ref{flat}), or by employing our assumptions on the \hge.

\ble\label{fflemma}
Under the assumptions of Proposition~\ref{forms}, for any 
$\phi\in\mbox{\em Hom}_\rho(V,\hO(B)P)$ there exist colinear homomorphisms
$\tilde{\phi}^\mu\in\mbox{\em Hom}_\rho(V,P)$ such that 
$\phi(v)=\sum_\mu\hl_\mu\tilde{\phi}^\mu(v),\; \fa v\in V$.
\ele\bpf
By assumption, we have
\[
\llp(\dsr\cc\phi)\ot id\lrp(v\0\ot v\1)
=\llp((\phi\te id)\ci\rho)\ot id\lrp(v\0\ot v\1),\;\mbox{i.e.,}
\]\[
\sum_\mu\hl_\mu\phi^\mu(v\0)\0\ot\phi^\mu(v\0)\1\ot v\1
=\sum_\mu\hl_\mu\phi^\mu(v\0)\ot v\1\ot v\2\, .
\]
Taking advantage of the faithful flatness of $P$, Theorem~I in~\cite{s-hj90}
and (1.6) in~\cite{d-y85} (Remark~3.3 in~\cite{s-hj90}),
 we know that there exists a unital colinear map  $j:H\ra P$. 
Applying 
\[
m_{\Omega(P)}\ci\llp id\ot(j\ci m)\lrp\ci(id\ot S\ot\id),
\]
where $m_{\Omega(P)}$ and $m$ are appropriate multiplication maps,
to both sides of the above equality, we get
\[
\sum_\mu\hl_\mu\phi^\mu(v\0)\0 j\llp S(\phi^\mu(v\0)\1)v\1\lrp
=\sum_\mu\hl_\mu\phi^\mu(v\0)j\llp S(v\1) v\2\lrp ,
\]
Hence, by the unitality of $j$, we obtain
\[
\phi(v)=\sum_\mu\hl_\mu\phi^\mu(v\0)\0 j\llp S(\phi^\mu(v\0)\1)v\1\lrp.
\]
On the other hand, using the colinearity of $j$ it is straightforward to
verify that each of the maps 
$v\st{\tilde{\phi}^\mu}{\longmapsto}\phi^\mu(v\0)\0 j\llp S(\phi^\mu(v\0)\1)v\1\lrp$
is colinear.
\epf\ \\

The next step is to take advantage
of the existence of the translation map $H\st{\tau}{\ra}P\ot_BP$, 
$\tau(h):=\chi^{-1}(1\ot h)$ (see Definition~\ref{hgdef}), and define
an auxiliary isomorphism 
\[
f:\hO(B)P\ra\hO(B)\ot_BP,~~ 
f:=(m\te_B id)\ci(id\te\tau)\ci\hD_{\cal R}\, .
\]
From the definition of the translation map it follows that 
\[
f(\hl p)=\hl p\0\tau(p\1)=\hl p\0\chi^{-1}(1\ot p\1)=\hl\chi^{-1}(p\0\ot p\1)
= \hl\chi^{-1}(\chi(1\ot_B p))=\hl\ot_Bp\, .
\]
(Note that $f$ is the inverse of the multiplication map.)
Moreover, let $I$ be the restriction to 
$\mbox{Hom}_\rho(V,\hO(B)P)$
of the canonical isomorphism from 
$\mbox{Hom}(V,\hO(B)P)$ to $\hO(B)P\ot V^*$.
Then we have a well-defined map 
\[
\hat{\ell}:=(id\ot_BI^{-1})\ci(f\te id)\ci I\;
:\;\mbox{Hom}_\rho(V,\hO(B)P)\lra\hO(B)P\ot_B\mbox{Hom}_\rho(V,P), 
\]
\[
\hat{\ell}(\phi)=\llp(id\ot_BI^{-1})\ci(f\te id)\lrp
\left(\sum_i\sum_\mu\hl_\mu\tilde{\phi}^\mu(e_i)\ot e^i\right)=
\sum_\mu\hl_\mu\ot_B\tilde{\phi}^\mu,
\]
 where  $\{e_i\}$ is a basis of $V$, $\{e^i\}$ its dual, and (by the above lemma)
we choose $\tilde{\phi}^\mu\in\mbox{Hom}_\rho(V,P)$   such that
$\phi(v)=\sum_\mu\lambda_\mu\tilde{\phi}^\mu(v)$. It is straightforward to check 
that $\hat{\ell}=\check{\ell}^{-1}$, as desired.
\epf

\bpr\label{flat}
Let $H$ be a Hopf algebra and $P\!\supseteq\! B$ an $H$\gal. Let $\check{\ell}$
be the map defined in Proposition~\ref{forms}. Then if $\hO(B)$ is flat as
a right $B$-module, $\check{\ell}$ 
is an isomorphism of graded left $\hO(B)$-modules.
\epr\bpf
Let $\check{\rho}: \mbox{Hom}(V,\hO(B)P)\lra\mbox{Hom}(V,\hO(B)P\ot H)$ be a
left $\hO(B)$-linear homomorphism defined by the formula
$\check{\rho}(\phi)(v)=\phi(v\0)\ot v\1-\phi(v)\0\ot\phi(v)\1$, and let 
$\tilde{\rho}$ denote its restriction to $\mbox{Hom}(V,P)$. Evidently, we have 
$\mbox{Ker}\,\check{\rho}=\mbox{Hom}_\rho(V,\hO(B)P)$ and 
$\mbox{Ker}\,\tilde{\rho}=\mbox{Hom}_\rho(V,P)$. Moreover, since $\hO(B)$ is
flat as a right $B$-module, we have the following commutative diagram with
exact rows of left $\hO(B)$-modules:
\beq\label{diag}
\def\normalbaselines{\baselineskip30pt
\lineskip3pt \lineskiplimit3pt }
\def\mapright#1{\smash{
\mathop{\!\!\!-\!\!\!\longrightarrow\!\!\!}
\limits^{#1}}}
\def\mapdown#1{\Big\downarrow
\rlap{$\vcenter{\hbox{$\scriptstyle#1$}}$}}
\matrix{
\! 0\!&\mapright{}&\!\hO(B)\te_B\mbox{Hom}_\rho(V,P)\!
&\mapright{}&\!\hO(B)\te_B\mbox{Hom}(V,P)\!
&\mapright{id\otimes_B\tilde{\rho}}&
\!\hO(B)\te_B\mbox{Hom}(V,P\te H)
\!\!\!\!\!\!\!\!\!\!\!\!\!\!\cr
&&\mapdown{\check{\ell}}&&\mapdown{\ell}&&\mapdown{\tilde{\ell}}\cr
\! 0\!&\mapright{}&\!\mbox{Hom}_\rho(V,\hO(B)P)\!
&\mapright{}&\!\mbox{Hom}(V,\hO(B)P)\!
&\mapright{\check{\rho}}&\!\mbox{Hom}(V,\hO(B)P\te H).
\!\!\!\!\!\!\!\!\!\!\!\!\!\!\cr
}
\eeq\ \\
Here $\ell$ is defined by the formula $\ell(\sum_\mu\hl_\mu\ot_B\phi^\mu)(v)
=\sum_\mu\hl_\mu\phi^\mu(v)$, and $\tilde{\ell}$ is given the same way.
With the help of the translation map $H\st{\tau}{\ra}P\ot_BP$, 
reasoning as in the proof of the preceding proposition, one can show that
$\ell$ and $\tilde{\ell}$ are isomorphisms. 
By standard diagram chasing (or completing the
left hand side of (\ref{diag}) with zeros and invoking the Five Isomorphism
Lemma), one can conclude from the diagram (\ref{diag}) that $\check{\ell}$
is also an isomorphism. 
\epf\ \\

If \ho\ is a strong connection form, then 
$
(id-\Pi^\omega)\ci d\ci\mbox{Hom}_\rho(V,P)\subseteq 
\mbox{Hom}_\rho(V,\hO^1(B)P).
$
Assuming also that 
the conditions allowing us to utilise one of the above propositions
are fulfilled, we can
define the {\em covariant derivative} associated to \ho\ in the
following way:
\beq\label{covde}
\nabla^\omega:\mbox{Hom}_\rho(V,P)\lra
\hO^{1}(B)\ot_B\mbox{Hom}_\rho(V,P),
~~
\nabla^\omega\xi:=
\check{\ell}^{-1}\llp(id-\Pi^\omega)\ci d\ci\xi\lrp . 
\eeq
One can check that $\nabla^\omega$ satisfies the Leibniz rule
$\nabla^\omega(b\xi)=b\nabla^\omega\xi+db\ot_B\xi$. Hence $\nabla^\omega$ can be
extended (by the Leibniz rule) to an endomorphism of 
$\hO(B)\ot_B\mbox{Hom}_\rho(V,P)$ which is of degree 1 with respect to the grading 
of~$\hO(B)$.

Our second group of results concerns the canonical connection on a
quantum principal homogeneous space (principal homogenous $H$\gal), 
which is the general
construction behind the Dirac $q$-monopole. A  principal homogenous
$H$\gal\ $B\inc P$ is a \hge\ obtained from a surjective Hopf algebra map 
$\pi:P\ra H$ which defines the right comodule structure by the formula 
$\dr:=(id\ot\pi)\ci\hD$. 
We know from the proof of \cite[Proposition~5.3]{bm93} that if
$B\inc P$ is a principal homogenous $H$\gal, and $i:H\ra P$
is a linear unital map such that $\pi\ci i=id$ (splitting of $\pi$)
and 
\beq\label{conad} 
\quad (id\tens\pi)\ci ad_R\circ i=(i\tens id)\ci ad_R\, ,
\eeq
then $\omega:=(S*d)\ci i$ is a connection form in the sense of 
Definition~\ref{confordef}. (Note that since $i$ is a splitting of a Hopf
algebra map, it is counital: $\he_H=\he_H\ci\pi\ci i=\he_P\ci i$.)
We call the thus constructed connection the {\em canonical connection (form)}
associated to splitting~$i$. (In what follows, 
we skip writing ``form" for the sake of brevity.)
Next step is towards a left-right symmetric
characterization of strong canonical connections.
\bpr\label{strocan}
The canonical connection associated to splitting
$i:H\ra P$ satisfying the above conditions
is strong {\em if and only if}
the splitting $i$ obeys in addition the right covariance condition
\[
(i\tens id)\ci\Delta = \Delta_R\ci i\, .
\]
\epr
\bpf 
First we need to reduce the strongness condition for the canonical
connection to a simpler form:
\ble\label{simfor}
The canonical connection $\ho$ associated to $i:H\ra P$ is strong 
{\em if and only if}
\beq\label{simeq}  
i(h\t)\t\tens h\o S\pi(i(h\t)\o)=  i(h)\tens 1,\quad \forall h\in H.
\eeq
\ele
\bpf
To simplify the notation, let us put $\pi(p)=\overline{p}$. Also, let $\Pi^\omega$ 
denote the connection associated to $\ho$, i.e., 
$\Pi^\omega(dp)=p\1\ho(\overline{p\2})\,$. (We take advantage of the fact that
$\dr=(id\ot\pi)\ci\hD$, see~(\ref{omefor}).)
Using the Leibniz rule we obtain:
\bea
&&
(id-\Pi^\omega)(dp)
\nonumber\\ &&
=d\,\llp p\1\, S(i(\overline{p\2})\1)\, i(\overline{p\3})\2\lrp
-p\1\, S(i(\overline{p\2})\1)\, d(i(\overline{p\3})\2)
\nonumber\\ &&
=d\,\llp p\1\, S(i(\overline{p\2})\1)\lrp\, i(\overline{p\3})\2
\nonumber\\ &&
=1\ot p-p\1\, S(i(\overline{p\2})\1)\ot i(\overline{p\3})\2
\, .\nonumber
\eea
On the other hand, applying $\dr\ot id$ to 
$p\1\, S(i(\overline{p\2})\1)\ot i(\overline{p\3})\2$
yields 
\[
p\o S(i(\overline{p\th})\t)\tens 
\overline{p\t}\;S\,\overline{i(\overline{p\th})\o}
\tens i(\overline{p\th})\th\, .
\]
Remembering that $(\hO^1B)P\inc B\ot P$, we conclude that
 the strongness condition (see Definition~\ref{scdef}, 
cf.~\cite[(11)]{m-s97}) of the canonical
connection is equivalent to
\[ 
p\o S(i(\overline{p\th})\t)\tens 
\overline{p\t}\;S\,\overline{i(\overline{p\th})\o}
\tens i(\overline{p\th})\th
=p\o S(i(\overline{p\t})\o)\tens 1\tens i(\overline{p\t})\t\, .
\]
The above equation is of the form $(id *f_1)(p)=(id *f_2)(p)$. 
Since the antipode
$S$ is the convolution inverse of $id$, it is equivalent to $f_1(p)=f_2(p)$. 
Therefore we can cancel the $p\o$ product from both sides. 
 Also, since $\pi$ is surjective and a
coalgebra map, we can replace $\pi(p)$ by a general element $h\in
H$.  Thus we arrive at
\[
S(i(h\2)\t)\tens h\1\;S\,\overline{i(h\2)\o}\tens i(h\2)\th
=S(i(h)\o)\tens 1\tens i(h)\t\, .
\]
Moreover,  for any Hopf algebra the map
$(S\tens id)\circ\Delta$ is injective (apply $\eps\tens id$). Consequently,
the strongness is equivalent to the condition
\[
i(h\2)\2\tens h\1\;S\,\overline{i(h\2)\1}
=i(h)\tens 1\, ,\;\;\;\forall\; h\in H\, ,
\]
as claimed.
\epf\ \\ \ \\
Note now that we can write the adjoint covariance of $i$, in an explicit manner, as
\beq\label{eman} 
i(h)\2\tens
S(\overline{i(h)\1})\,\overline{i(h)\3}=i(h\2)\tens (Sh\1) h\3
\, ,\quad\forall h\in H.
\eeq 
In this case
\bea
&&
i(h\1)\tens h\2
\nonumber\\ &&
=i(h\3)\tens h\1 S(h\2)h\4
\nonumber\\ &&
=(1\ot h\1)((i\ot id)\cc ad_R)(h\2)
\nonumber\\ &&
=i(h\2)\t\tens
h\1 S(\overline{i(h\2)\1})\,\overline{i(h\2)\3}
\, .\nonumber
\eea 
Assume that $\ho$ is strong. Hence, by the above lemma, the strongness condition
implies that  $i(h\o)\tens h\t=i(h)\o\tens\overline{i(h)\t}$ as
required. 
Conversely, using the right covariance of $i$ for the first step and (\ref{eman}) 
for the second, we compute the left hand
side of (\ref{simeq}) as
\align{ \equad && i(h\t)\t\tens h\o S(\overline{i(h\t)\1})
 \,\overline{i(h\t)\th}\, S(\overline{i(h\t)\4})\\
&&=i(h\t)\t\tens h\o\, S(\overline{i(h\t)\o})\,\overline{i(h\t)\th}\, S h\th\\ &&=
i(h\3)\tens h\o S(h\t) h\4 S(h\5)\\ &&=i(h)\tens 1.} 
Hence the canonical connection is strong by Lemma~\ref{simfor}.
\epf 
\bco\label{equiv} 
Assume that antipode $S$ is injective. Then
strong canonical connections are in 1-1 correspondence with
linear unital splittings of $\pi$ obeying the two conditions
\[(i\tens id)\ci\Delta = \Delta_R\ci i,\quad (id\tens i)\ci\Delta =
\Delta_L\ci i,\]
where $\Delta_R=(id\tens\pi)\ci\Delta,\;
\Delta_L=(\pi\tens id)\ci\Delta$.
\eco
\bpf 
Assume first that the canonical connection associated to $i$ is strong. Then,
by the preceding proposition, $i$ is right covariant and (\ref{eman}) holds. Hence
\align{\equad && i(h\o)\t\tens S(\overline{i(h\o)\o})\, h\t\\ &&=i(h)\o\t\tens
S(\overline{i(h)\o\o})\,\overline{i(h)\t}\\ &&=i(h)\t\tens S(\overline{i(h)\o})\,
\overline{i(h)\th}\\ &&=i(h\t)\tens (Sh\o)h\th.} 
Reasoning as in the proof of Lemma~\ref{simfor}, we can cancel $h\t$ and $h\th$ from
the two sides. Then cancelling $S$ from both sides (we assume $S$ to be injective), 
we have
$i(h)\t\tens\overline{i(h)\o}=i(h\t)\tens h\o$, which is the left
covariance condition. 

Conversely, if the left and right covariance
conditions hold then
\bea
&&
i(h\t)\tens
(Sh\o)h\th\equad 
\nonumber\\ &&
=i(h\t\o)\tens (Sh\o)h\t\t
\nonumber\\ &&
=i(h\t)\o\tens
(Sh\o)\overline{i(h\t)\t}
\nonumber\\ &&
=i(h)\t\o\tens S(\overline{i(h)\o})\,\overline{i(h)\t\t}
\, ,\nonumber
\eea
which is the same as (\ref{eman}). Invoking again the preceding proposition,
we can conclude that the canonical connection associated to $i$ is strong
as required.
\epf 

\bre\label{counit}\em
Let $\pi:P\ra H$ be a Hopf algebra surjection. If a linear map $i:H\ra P$ is
counital and left or right colinear, then $i$ is a splitting of $\pi$, i.e.,
$\pi\ci i=id$. Indeed, if $i$ is right colinear ($i(h)\1\ot\pi(i(h)\2)=
i(h\1)\ot h\2$), we have:
\[
(\pi\ci i)(h)=
\he\llp(\pi\ci i)(h)\1\lrp(\pi\ci i)(h)\2=
\he\llp\pi(i(h)\1)\lrp\pi(i(h)\2)=
\he\llp\pi(i(h\1))\lrp h\2=
h.
\]
The left-sided case is analogous.
\ere

We end this section by showing how to obtain a projector matrix (explicit embedding
of a projective module in a free module) from the canonical strong connection.
It is known~\cite{dh98} that strong connection forms on $P$ are {\em equivalent}
to unital left $B$-linear right $H$-colinear splittings of the multiplication map
$m: B\ot P\ra P$. Explicitly, if \ho\ is a strong connection form, then
\beq\label{s}
s:P\lra B\ot P,~~~ s(p)=p\ot 1 +p\0\ho(p\1) 
\eeq
gives the desired splitting. (Solving this equation for \ho\ one gets
$\ho(h)=h^{[1]}s(h^{[2]})-1\ot\he(h)$, where 
$h^{[1]}\ot_B h^{[2]}=\chi^{-1}(1\ot h)$, summation understood, see 
Definition~\ref{hgdef}.) 
In particular, for the canonical
strong connection associated to a bicovariant splitting $i$ 
(i.e., $\ho=(S*d)\ci i$), we have:
\beq\label{i}
s(p)=p\1 Si(\overline{p\2})\1\ot i(\overline{p\2})\2\, .
\eeq 
Note that a splitting of the multiplication map is almost the same as a projector
matrix, for it is an embedding of $P$ in the free $B$-module $B\ot P$. (We will
use formula (\ref{i}) in the next section to compute projector matrices of quantum 
Hopf line bundles from the Dirac $q$-monopole connection.) To turn (\ref{i}) into
a concrete recipe for producing finite size projector matrices of finitely
generated projective modules, let us claim the following general lemma:
\ble\label{gen}
Let $A$ be an algebra and $M$ a projective left $A$-module generated by linearly
independent generators $g_1,...,g_n$. Also, let $\{\widetilde{g}_\mu\}_{\mu\in I}$
be a completion of $\{g_1,...,g_n\}$ to a linear basis of $M$, $f_2$ be a left
$A$-linear splitting of the multiplication map $m:A\ot M\ra M$ given by the 
formula 
$f_2(g_k)=\sum_{l=1}^{n}a_{kl}\ot g_l+\sum_{\mu\in I}a_{k\mu}\ot\widetilde{g}_\mu$,
and  $c_{\mu l}\in A$ a choice of coefficients such that 
$\widetilde{g}_\mu=\sum_{l=1}^{n}c_{\mu l}g_l$. Then 
$e_{kl}=a_{kl} +\sum_{\mu\in I}a_{k\mu}c_{\mu l}$ defines a projector matrix
of $M$, i.e., $e\in M_n(A)$, $e^2=e$ and $A^ne$ and $M$ are isomorphic as left 
$A$-modules.
\ele\bpf
Note first that we do not lose any generality by assuming $g_1,...,g_n$ to be
linearly independent (we can always remove generators that are linear combinations
of other generators), and that a splitting of the multiplication map always
exists by the projectivity assumption (cf.~\cite[Section~8]{cq95}). Let $N$ be
the kernel of the surjection $f_1:A^n\ra M=A^n/N$, $f_1(e_k)=g_k$,
$k\in\{1,...,n\}$, where $\{e_k\}_{k\in\{1,...,n\}}$ is the standard basis of $A^n$,
i.e., $e_k$ is the row with zeros everywhere except for the $k$-th place where
there is 1.
We have the following commutative diagram of left $A$-module homomorphisms whose
rows are exact:
\vspace*{3mm}\beq\label{gend}
\def\normalbaselines{\baselineskip30pt\lineskip3pt\lineskiplimit3pt}
\def\mapright#1{\smash{\mathop{-\!\!\!\lra}\limits^{#1}}}
\def\Mapright#1{\smash{\mathop{-\!\!\!-\!\!\!-\!\!\!-\!\!\!\lra}\limits^{#1}}}
\def\mapleft#1{\smash{\mathop{\longleftarrow\!\!\!-}\limits_{#1}}}
\def\mapdown#1{\Big\downarrow\rlap{$\vcenter{\hbox{$\scriptstyle#1$}}$}}
\def\mapup#1{\Big\uparrow\llap{$\vcenter{\hbox{$\scriptstyle#1$}}$}}
\matrix{
0 &\mapright{}& A\ot N &\mapright{}& A\ot A^n 
&\stackrel{\large\Mapright{id\ot f_1}}{\mapleft{f_3}}& 
A\ot M &\mapright{}& 0 \cr
&& && \mapdown{f_4} && \mapup{f_2~~}\mapdown{m} && \cr
0 &\mapright{}& N &\mapright{}& A^n &\mapright{f_1}& M &\mapright{}& 0 \cr
}
\eeq
Here $f_2$ is a splitting of the multiplication map ($m\ci f_2=id$), $f_3$
a splitting of $id\ot f_1$ (which exists because $A\ot M$ is free), and $f_4$
is the multiplication map on $A\ot A^n$. From the commutativity of the diagram
we can infer that $f_4\ci f_3\ci f_2$ is a splitting of $f_1$:
\[
f_1\ci f_4\ci f_3\ci f_2=m\ci(id\ot f_1)\ci f_3\ci f_2=id.
\]
Hence $f_e:=f_4\ci f_3\ci f_2\ci f_1$ is an idempotent ($f_e^2=f_e$) and $f_e(A^n)$
is isomorphic to~$M$, as needed. To compute a matrix of $f_e$,
we choose a splitting $f_3$ so that 
$f_3(1\ot g_k)=1\ot e_k$, 
$f_3(1\ot \widetilde{g}_\mu)=1\ot \sum_{l=1}^{n}c_{\mu l}e_l$, 
$\sum_{l=1}^{n}c_{\mu l}g_l=\widetilde{g}_\mu$, $k\in\{1,...,n\}$, $\mu\in I$.
Then
\begin{eqnarray*}\label{gena}
f_e(e_k)\!\!\!\!\!\!&&=(f_4\ci f_3\ci f_2)(g_k)
\\ && =
(f_4\ci f_3)
\llp\sum_{l=1}^{n}a_{kl}\ot g_l+\sum_{\mu\in I}a_{k\mu}\ot\widetilde{g}_\mu\lrp
\\ && =
f_4\llp
\sum_{l=1}^{n}a_{kl}\ot e_l+\sum_{\mu\in I}a_{k\mu}\ot\sum_{l=1}^{n}c_{\mu l}e_l
\lrp
\\ && =
\sum_{l=1}^{n}\llp a_{kl}+\sum_{\mu\in I}a_{k\mu}c_{\mu l}\lrp e_l\, .
\end{eqnarray*}
This means that $(a_{kl}+\sum_{\mu\in I}a_{k\mu}c_{\mu l})_{k,l\in\{1,...,n\}}$ 
is a projector matrix of $M$, as claimed.
\epf\\ ~\\
Observe that if $a_{k\mu}=0$ for all $k$ and $\mu$, the matrix elements of $e$
are simply $a_{kl}$, and can be directly read off from the formula for splitting
$f_2$ written in terms of the module generators $g_1,...,g_n$. By a completely
analogous reasoning, the same kind of lemma is true for right modules.

\section{Projective module form of the Dirac q-monopole}

Recall that $A(SL_q(2))$ is a Hopf algebra over a field $k$
generated by $1,\, \ha,\, \hb,\, \hg,\, \hd$, satisfying the 
following relations: 
\bea \label{comm}
&& \ha\hb=q^{-1} \hb\ha~,~~\ha\hg=q^{-1}\hg\ha~,~~\hb\hd=q^{-1}
\hd\hb~,~~\hb\hg=\hg\hb~,
~~\hg\hd=q^{-1}\hd\hg~, \nonumber \\ 
&& \ha\hd-\hd\ha=(q^{-1}-q)\hb\hg~,~~\ha\hd-q^{-1}\hb\hg=\hd\ha-q\hb\hg=1\, ,  
\eea
where $q\in k\setminus\{0\}$. 
The comultiplication $\Delta$, counit $\varepsilon$, and antipode
$S$ of \slq\ are defined by the following formulas:
\[ 
\Delta\pmatrix{\alpha & \beta \cr \gamma & \delta \cr}=
\pmatrix{\alpha\ot1 & \beta\ot1 \cr \gamma\ot1 & \delta\ot1 \cr} 
\pmatrix{1\ot\alpha & 1\ot\beta \cr 1\ot\gamma & 1\ot\delta \cr},
\]\[
\varepsilon\pmatrix{\alpha & \beta \cr \gamma & \delta\cr}=
\pmatrix{ 1 & 0 \cr 0 & 1 \cr},~~~
S\pmatrix{\alpha & \beta \cr \gamma & \delta \cr}=
\pmatrix{\delta &-q\beta \cr -q^{-1}\gamma & \alpha \cr}.
\]

Now we need to recall the construction of the standard quantum sphere of Podle\'s
and the quantum principal Hopf fibration.
The standard quantum sphere is singled out among the
principal series of Podle\'s quantum spheres by the property that it can be
constructed as a quantum quotient space~\cite{p-p87}. In algebraic terms it
means that its coordinate ring can be obtained as the subalgebra of coinvariants
of a comodule algebra. To carry out this construction, first we need the
right coaction on \slq\ of the commutative and cocommutative Hopf algebra 
$k[z,z^{-1}]$ generated by the grouplike element $z$ and its inverse. This
Hopf algebra can be obtained as the quotient of \slq\ by the Hopf ideal generated
by the off-diagonal generators $\hb$ and~$\hg$. Identifying the image of \ha\ and
\hd\ under the Hopf algebra surjection $\pi:\slq\ra k[z,z^{-1}]$ with $z$ and 
$z^{-1}$ respectively, we can describe the right coaction $\dr:=(id\ot\pi)\ci\hD$
by the formula:
\[
\dr\pmatrix{\ha & \hb \cr \hg & \hd \cr}=
\pmatrix{\ha\ot z & \hb\ot z^{-1} \cr \hg\ot z & \hd\ot z^{-1} \cr}.
\]
We call the subalgebra of coinvariants defined by this coaction the coordinate
ring of the (standard) {\em quantum sphere}, and denote it by \asq. Since 
\[
k[z,z^{-1}]=\slq/(\asq\cap\mbox{Ker}\he)\slq
\] 
by Remark~\ref{local}, we know from the general argument that $\asq\inc\slq$ is
a principal homogenous $k[z,z^{-1}]$\gal. (If $P$ is a Hopf algebra, $I$ a Hopf
ideal, $B$ the subalgebra of coinvariants under the coaction $\dr=(id\ot\pi)\ci\hD$,
$P\stackrel{\pi}{\ra}P/I$, and $I=(B\cap\mbox{Ker}\he)P$, then we can define the 
inverse of the canonical map
by $\chi^{-1}(p'\ot\pi(p))=p'Sp\1\ot_Bp\2$.) We refer to the quantum principal
bundle given by this \hge\ as the
{\em quantum principal Hopf fibration}. (An $SO_q(3)$ version of this quantum 
fibration was studied in~\cite{bm93}.)

The main point of this section is to compute projector matrices of quantum
Hopf line bundles associated to the just described Hopf $q$-fibration.
\bde\label{line}
Let $\rho_n:k[z,z^{-1}]\ra k\ot k[z,z^{-1}]$, $\rho_n(1)=1\ot z^{-n}$, 
$n\in{\Bbb Z}$, be a one-dimensional corepresentation of $k[z,z^{-1}]$.
We call the \asq-bimodule of colinear maps $\mbox{\em Hom}_{\rho_n}(k,\slq)$
the (bimodule of) {\em quantum Hopf line bundle of winding number $n$}.
\ede
Since we deal here with one-dimensional corepresentations, we identify colinear
maps with their value at~1. We have
\[
\mbox{Hom}_{\rho_n}(k,\slq)\tilde{=}\{p\in\slq\ |\;\dr p=p\ot z^{-n}\}=:P_n
\]
as \asq-bimodules. With the help of the PBW basis 
$\ha^k\hb^l\hg^m,\;\hb^p\hg^r\hd^s,\; k,l,m,p,r,s\in\IN_0,\; k>0$ of \slq,
one can show that 
\begin{eqnarray*}\label{osum}
&&
P_n=\left\{
\begin{array}{ll}
\sum_{k=0}^{-n}\asq\ \ha^{-n-k}\hg^{k}=\sum_{k=0}^{-n}\ha^{-n-k}\hg^{k}\asq
& \mbox{for $n\leq 0$}\\
\sum_{k=0}^{n}\asq\ \hb^{k}\hd^{n-k}=\sum_{k=0}^{n}\hb^{k}\hd^{n-k}\asq
& \mbox{for $n\geq 0$}, 
\end{array}
\right.
\end{eqnarray*}
and $\slq=\bigoplus_{n\in{\Bbb Z}}P_n$ (cf.~\cite[(1.10)]{mmnnu91}).

Next, similarly to \cite{bm93}, we consider the canonical connection
induced by the bicovariant splitting $i(z^n)=\alpha^n,\ i(z^{-n})=\delta^n$
(see~\cite{bm97}). By Corollary~\ref{equiv} it induces a
strong connection. We call this connection the (Dirac) {\em $q$-monopole}.
Now, formula (\ref{i}) gives us a splitting $s:\slq\ra\asq\ot\slq$, and we
can claim:
\bpr\label{ml}
Put 
\begin{eqnarray*}\label{ekl}
&&
(e_n)_{kl}=\left\{
\begin{array}{ll}
\ha^{-n-k}\hg^{k}{\scriptsize\pmatrix{-n\cr l\cr}}_{\!q^2}(-q)^l\hb^l\hd^{-n-l}
& \mbox{for $n\leq 0$}\\
\hb^{k}\hd^{n-k}{\scriptsize\pmatrix{n\cr l\cr}}_{\!q^2}(-q)^{-l}\ha^{n-l}\hg^l
& \mbox{for $n\geq 0$}. 
\end{array}
\right.
\end{eqnarray*}
Then, for any $n\in\IZ$, $e_n\in M_{|n|+1}(\asq)$, $e_n^2=e_n$, and 
$\asq^{|n|+1}e_n$ is isomorphic to $P_n$ as a left \asq-module. 
\epr\bpf
Recall first that if $qxy=yx$, then 
$(x+y)^n=\sum_{k=0}^n{\scriptsize\pmatrix{n\cr k\cr}}_{\!q}x^ky^{n-k}$, where
\[
{\scriptsize\pmatrix{n\cr k\cr}}_{\!q}=
\frac{(q-1)...(q^n-1)}{(q-1)...(q^k-1)(q-1)...(q^{n-k}-1)}
\]
 are the $q$-binomial coefficients. (See, e.g.,~\cite[p.85]{m-s95}.)
Taking advantage
of formula  (\ref{i}) in the $q$-monopole case, we compute:
\begin{eqnarray*}
s(\ha^{m-k}\hg^k)\!\!\!&&=\ha^{m-k}\hg^kSi(z^m)\1\ot i(z^m)\2
\\ && = 
\sum_{l=0}^m\ha^{m-k}\hg^{k}{\scriptsize\pmatrix{m\cr l\cr}}_{\!q^2}
S(\ha^{m-l}\hb^{l})\ot \ha^{m-l}\hg^{l}
\\ && = 
\sum_{l=0}^m\ha^{m-k}\hg^{k}{\scriptsize\pmatrix{m\cr l\cr}}_{\!q^2}(-q)^l
\hb^{l}\hd^{m-l}\ot \ha^{m-l}\hg^{l}.
\end{eqnarray*}
Similarly, $s(\hb^k\hd^{n-k})=\sum_{l=0}^n
\hb^{k}\hd^{n-k}{\scriptsize\pmatrix{n\cr l\cr}}_{\!q^2}(-q)^{-l}\ha^{n-l}\hg^l
\ot \hb^{l}\hd^{n-l}$. Thus we have verified that $s$ preserves the direct sum
decomposition of \slq, i.e., $s(P_n)\inc\asq \ot P_n$, $n\in\IZ$. Hence, 
by restriction, we have a splitting of the left multiplication map for each~$P_n$.
The claim of the proposition follows directly from Lemma~\ref{gen} and the above
formulas for $s$.
\epf
\bre\label{direct}\em
Observe that for $n\geq 0$ we can write $e_n=uv^T$, where 
$u^T=(\hd^n,...,\hb^k\hd^{n-k},...,\hb^n)$ and $v^T=(S(\hd^n),...,
{\scriptsize\pmatrix{n\cr k\cr}}_{\!q^2}S(\hg^k\hd^{n-k}),...,S(\hg^n))$.
Since 
\[
v^Tu
=\sum_{k=0}^{n}{\scriptsize\pmatrix{n\cr k\cr}}_{\!q^2}S(\hg^k\hd^{n-k})\hb^k\hd^{n-k}
=S((\hd^n)\1)(\hd^n)\2=\he(\hd^n)=1,
\]
we can directly see that $e_n^2=e_n$. The case $n\leq 0$ is similar.
\ere
\bre\label{local}\em
We can define the fibre
of a quantum vector bundle over a classical point (understood as a number-valued
 algebra homomorphism) as the localization of the
module of ``sections" of this bundle at the kernel of this homomorphism. 
The standard Podle\'s quantum sphere that we consider  here has one classical
point given by the restriction of the counit map~\he. Let us consider the
quantum Hopf line bundles as left \asq-modules~$P_n$. We can then regard
the localization $P_n/\asq^+P_n$, $\asq^+:=\mbox{Ker}\he\cap\asq$, as the fibre vector
space of $P_n$ over the point given by~$\asq^+$. (Note that $P_n/\asq^+P_n$ is 
automatically a vector space over $\asq/\asq^+=k$.) Since $\he(\asq^+P_n)=0$, \he\
induces a linear map $\tilde{\he}:P_n/\asq^+P_n\ra k$ given by the formula
$\tilde{\he}(p/\asq^+P_n)=\he(p)$. Assume now that $n\geq 0$. Arbitrary $p\in P_n$
can be written as $p=\sum_{l=0}^{n}b_l\hb^l\hd^{n-l}$, $b_l\in\asq$. Hence
$\tilde{\he}(p/\asq^+P_n)=\he(b_0)$, and we can conclude that $\tilde{\he}$ is 
surjective. Note now that $\hb=(-q^{-1}\hb\hg)\hb+(q\ha\hb)\hd$, and consequently,
for $l>0$, 
$\hb^l\hd^{n-l}=(-q^{-1}\hb\hg)\hb^l\hd^{n-l}+(q\ha\hb)\hd\hb^{l-1}\hd^{n-l}
\in\asq^+P_n$. It follows that 
\begin{eqnarray*}
(\sum_{l=0}^{n}b_l\hb^l\hd^{n-l})/\asq^+P_n
\!\!\!\!\!\! &&
=b_0\hd^{n}/\asq^+P_n
\\ &&
=\he(b_0)\hd^{n}/\asq^+P_n+(b_0-\he(b_0))\hd^{n}/\asq^+P_n
\\ &&
=\he(b_0)\hd^{n}/\asq^+P_n.
\end{eqnarray*}
This entails the injectivity of $\tilde{\he}$. Thus $\tilde{\he}$ is an isomorphism,
 and we can infer that the fibre $P_n/\asq^+P_n$ is a one-dimensional vector space,
exactly as expected for a line bundle. The reasoning for $n\leq 0$ is analogous,
and relies on the identity $\hg=(-q\hb\hg)\hg+(q^{-1}\hd\hg)\ha$. 
This agrees with the fact that $\slq=\bigoplus_{n\in\IZ}P_n$ and
$\slq/\asq^+\slq=k[z,z^{-1}]=\bigoplus_{n\in\IZ}kz^n$. The latter equality can be
directly seen as follows: Since \hb\ and \hg\ $q$-commute with all monomials,
the two-sided ideal $\langle\hb,\hg\rangle=\hb\slq+\hg\slq$. Thus, as 
$\hb,\hg\in\asq^+\slq$ by the above formulas, we have 
$\langle\hb,\hg\rangle\inc\asq^+\slq$. On the other hand, since $\asq^+$ is the ideal
in \asq\ generated by \ha\hb, \hb\hg, \hg\hd, we also have 
$\asq^+\slq\inc\langle\hb,\hg\rangle$. Hence 
$k[z,z^{-1}]=\slq/\langle\hb,\hg\rangle=\slq/\asq^+\slq$.
\ere

To compute projector matrices of the quantum Hopf line bundles thought of
as right \asq-modules, we need a right-sided version of formula~(\ref{i}).
A natural first candidate appears to be:
\beq\label{ii}
\tilde{s}(p)=i(\overline{p\1})\1\ot S(i(\overline{p\1})\2)p\2\, .
\eeq
It is evidently a splitting of the multiplication map $m:\slq\ot\slq\ra\slq$.
Only now it is right linear under left coinvariants. By left coinvariants we
understand here $\tasq\ :=\{p\in \slq\ |\; \hD_Lp=1\ot p\}$, where 
$\hD_L=(\pi\ot id)\ci\hD$. On generators, we have explicitly:
\[
\hD_L\pmatrix{\ha & \hb \cr \hg & \hd \cr}=
\pmatrix{z\ot\ha & z\ot\hb \cr z^{-1}\ot\hg & z^{-1}\ot\hd \cr}.
\]
Using the PBW basis 
$\ha^k\hb^l\hg^m,\;\hb^p\hg^r\hd^s,\; k,l,m,p,r,s\in\IN_0,\; k>0$ of \slq,
one can show that \tasq\ is a unital subalgebra of \slq\ generated by
\ha\hg, \hb\hd, \hb\hg. We want to prove now that the image of $\tilde{s}$
lies in $\slq\ot\tasq$. To this end we note that the right covariance of $i$
implies the formula 
$
i(h)\1\ot\overline{i(h)\3}\ot i(h)\2=i(h\1)\1\ot h\2\ot i(h\1)\2\, .
$
With the above formula at hand, one can verify that
$
((id\ot\hD_L)\ci\tilde{s})(p)=i(\overline{p\1})\1\ot 1\ot 
S(i(\overline{p\1})\2)p\2\, ,
$
as needed. Thus we can conclude that $\tilde{s}$ is a right \tasq-linear
splitting of the multiplication map $\slq\ot\tasq\ra\slq$. However, \tasq\
and \asq\ are different subalgebras of \slq, and we want to find projector matrices
for $P_n$ thought of as right \asq-modules. To our aid comes the transpose
automorphism of \slq\ defined on generators by
\[
T\pmatrix{\ha & \hb \cr \hg & \hd \cr}=
\pmatrix{\ha & \hg \cr \hb & \hd \cr}.
\]
One can check directly that $T$ is well defined. In particular, when we work over 
$\Bbb C$, \slq\ has a natural $*$-algebra structure for $q$ real, namely
\[ 
*\pmatrix{\ha & \hb \cr \hg & \hd \cr}=
\pmatrix{\hd & -q^{-1}\hg \cr -q\hb & \ha \cr},
\]
and we can simply define $T=*\ci S$. This automorphism gives an isomorphism
between \asq\ and \tasq. We have $T(\asq)=\tasq$ and $T(\tasq)=\asq$. (Note
that $T^2=id$.) It is straightforward to verify that 
$\check{s}:=(T\ot T)\ci\tilde{s}\ci T$ is a right \asq-linear splitting of
the right multiplication map $m:\slq\ot\asq\ra\slq$. We can now proceed as in
the left-sided case to prove:
\bpr\label{mr}
Put 
\begin{eqnarray*}\label{elk}
&&
(f_n)_{lk}=\left\{
\begin{array}{ll}
{\scriptsize\pmatrix{-n\cr l\cr}}_{\!q^2}(-q)^{-l}\hb^l\hd^{-n-l}\ha^{-n-k}\hg^{k}
& \mbox{for $n\leq 0$}\\
{\scriptsize\pmatrix{n\cr l\cr}}_{\!q^2}(-q)^{l}\ha^{n-l}\hg^l\hb^{k}\hd^{n-k}
& \mbox{for $n\geq 0$}. 
\end{array}
\right.
\end{eqnarray*}
Then, for any $n\in\IZ$, $f_n\in M_{|n|+1}(\asq)$, $f_n^2=f_n$, and 
$f_n\asq^{|n|+1}$ is isomorphic to $P_n$ as a right \asq-module. 
\epr\bpf
We have:
\begin{eqnarray*}
\check{s}(\ha^{m-k}\hg^k)\!\!\!&&=(T\ot T)(\tilde{s}(\ha^{m-k}\hb^k))
\\ && =
(T\ot T)\llp i(z^m)\1\ot S(i(z^m)\2)\ha^{m-k}\hb^k \lrp
\\ && = 
(T\ot T)\llp\sum_{l=0}^m{\scriptsize\pmatrix{m\cr l\cr}}_{\!q^2}\ha^{m-l}\hb^{l}
\ot S(\ha^{m-l}\hg^{l})\ha^{m-k}\hb^{k}\lrp
\\ && = 
(T\ot T)\llp\sum_{l=0}^m\ha^{m-l}\hb^{l}\ot {\scriptsize\pmatrix{m\cr l\cr}}_{\!q^2}
(-q)^{-l}\hg^{l}\hd^{m-l}\ha^{m-k}\hb^{k}\lrp
\\ && = 
\sum_{l=0}^m\ha^{m-l}\hg^{l}\ot {\scriptsize\pmatrix{m\cr l\cr}}_{\!q^2}
(-q)^{-l}\hb^{l}\hd^{m-l}\ha^{m-k}\hg^{k}\, .
\end{eqnarray*}
Similarly, $\check{s}(\hb^k\hd^{n-k})=\sum_{l=0}^n
\hb^{l}\hd^{n-l}\ot {\scriptsize\pmatrix{n\cr l\cr}}_{\!q^2}(-q)^{l} 
\ha^{n-l}\hg^lb^{k}\hd^{n-k}$. Hence $\check{s}(P_n)\inc P_n\ot\asq$, $n\in\IZ$.  
By restriction of $\check{s}$, 
we have a splitting of the right multiplication map for each~$P_n$.
The claim of the proposition follows from the right-sided version of 
Lemma~\ref{gen} and the above formulas for~$\check{s}$.
\epf\\

Finally, let us observe that, identifying $\mbox{Hom}_{\rho_n}(k,\slq)$ with
$P_n$, we can view the covariant derivative 
$\nabla^\omega_n:\mbox{Hom}_{\rho_n}(k,\slq)\ra
\hO^1\asq\ot_{A(S^2_q)}\mbox{Hom}_{\rho_n}(k,\slq)$
associated to the $q$-monopole by~(\ref{covde}), as the Grassmannian connection 
associated to the splitting $s_n:=s|_{P_n}$. More precisely, let 
$
\psi:\mbox{Hom}_{\rho_n}(k,\slq)\ra P_n,\; \psi(\xi)=\xi(1)
$
be the identification isomorphism mentioned above. The Grassmannian connection
associated to the splitting $s_n:P_n\ra\asq\ot P_n$  is by 
definition the connection $\tilde{\nabla}_n^s:P_n\ra\hO^1A(S^2_q)\ot P_n$ given by the
formula $\tilde{\nabla}_n^sp=\sum_i db_i\ot_{A(S^2_q)} p_i$, where 
$\sum_i b_i\ot p_i:=s(p)$.
(See \cite[(54)]{cq95} or \cite[(8.27)]{l-g97}
 for the right-sided version.) We want to show that
\[
\nabla_n^\omega =(id\ot_{A(S^2_q)}\psi^{-1})\ci\tilde{\nabla}_n^s\ci\psi,
\; n\in\IZ,
\]
or equivalently that
\[
\Lall\xi\in\mbox{Hom}_{\rho_n}(k,\slq),\, n\in\IZ:\; 
\llp\check{\ell}(\nabla_n^\omega\xi)\lrp(1)=
\llp
(\check{\ell}\ci(id\ot_{A(S^2_q)}\psi^{-1})\ci\tilde{\nabla}_n^s\ci\psi)
(\xi)\lrp(1).
\]
(See Proposition~\ref{forms} and~(\ref{covde}).) 
Notice that we can use here either Proposition~\ref{forms} or Proposition~\ref{flat}
to guarantee that $\nabla^\omega_n,\; n\in\IZ$, makes sense.
Indeed, since $k[z,z^{-1}]$ admits the Haar functional ($h_H:k[z,z^{-1}]\ra k$,
$h_H(z^n)=\hd_{0n}$), we can construct a unital right colinear mapping 
$j:k[z,z^{-1}]\ra\slq$, $j:=\eta\ci h_H$, where $\eta: k\ra\slq$ is the unit map,
so that \slq\ is injective as a right $k[z,z^{-1}]$-comodule.
Thus, as the antipode of $k[z,z^{-1}]$ is bijective,
  \slq\ is left and right faithfully flat over \asq\ by \cite[Theorem~I]{s-hj90},
and Proposition~\ref{forms} applies. (In fact, we used the existence of 
a unital right colinear mapping to prove Proposition~\ref{forms}.)
Also, $\hO^1\asq$ is isomorphic with $\asq/k\ot\asq$ as a right \asq-module via
$db\to b/k\ot 1$, so that it is free, whence flat. Therefore Proposition~\ref{flat}
applies as well.
Now, we put $s(\xi(1))=b_i\ot\xi(1)_i$,
$\xi_i(1)=\xi(1)_i$, and taking advantage of $m\ci s_n=id$, (\ref{s}), (\ref{omefor}) 
and (\ref{covde}) compute:
\begin{eqnarray*}\label{grass}
&&
\llp
(\check{\ell}\ci(id\ot_{A(S^2_q)}\psi^{-1})\ci\tilde{\nabla}_n^s\ci\psi)
(\xi)\lrp(1)
\\ && =
\llp(\check{\ell}\ci(id\ot_{A(S^2_q)}\psi^{-1})\lrp
(\sum_i db_i\ot_{A(S^2_q)} \xi(1)_i)(1)
\\ && =
\sum_i (db_i)\xi(1)_i
\\ && =
1\ot (m\ci s_n)(\xi(1)) - s_n(\xi(1))
\\ && =
1\ot \xi(1) - \xi(1)\ot 1 - \xi(1)\0\ho(\xi(1)\1)
\\ && =
d\xi(1) - \Pi^\omega(d\xi(1))
\\ && =
\llp\check{\ell}(\nabla_n^\omega\xi)\lrp(1).
\end{eqnarray*}
This is exactly as one should expect, since we have constructed the splitting
$s:\slq\ra\asq\ot\slq$ from the connection form \ho\ by formula~(\ref{s}).

\section{Chern-Connes pairing for the \boldmath$n=-1$ bimodule}

The aim of this section is to compute the left and right Chern numbers of the
left and right finitely generated projective bimodule $P_{-1}$ 
describing the quantum Hopf line bundle of winding number~$-1$. 
This computation is a simple example 
of the Chern-Connes pairing between $K$-theory and cyclic 
cohomology~\cite{c-a94,l-jl97}. 

To obtain the desired Chern numbers we need 
to evaluate (to pair) the appropriate even cyclic cocycle with the left and right
projector matrix respectively. Since the positive even cyclic cohomology
$HC^{2n}(\asq)$, $n>0$, is the image of the periodicity operator applied to
$HC^{0}(\asq)$, and the pairing is compatible with the action of the periodicity
operator, the even cyclic cocycle computing Chern numbers is necessarily of
degree zero, i.e., a trace. This trace is explicitly provided 
in~\cite[(4.4)]{mnw91}. Adapting \cite[(4.4)]{mnw91} to our special case of the
standard Podle\'s quantum sphere, we obtain: 
\bea\label{tr}
&&
\tau^1((\ha\hb)^m\hz^n)=\left\{
\begin{array}{ll}
(1-q^{2n})^{-1} & \mbox{for $n>0$, $m=0$,}\\
0 & \mbox{otherwise}, 
\end{array}
\right.
\nonumber\\ && \ \nonumber\\ &&
\tau^1((\hg\hd)^m\hz^n)=\left\{
\begin{array}{ll}
(1-q^{2n})^{-1} & \mbox{for $n>0$, $m=0$,}\\
0 & \mbox{otherwise}, 
\end{array}
\right.
\eea
where $\hz:=-q^{-1}\hb\hg$.

The fact that the ``Chern cyclic cocycle" is in degree zero is a quantum effect
caused by the non-classical structure of $HC^*(\asq)$
(see~\cite{mnw91}). 
In the
classical case the corresponding cocycle is in degree two, as it comes from the
volume form of the two-sphere.

Since $\tau^1$ is a 0-cyclic cocycle, the pairing is given by the formula
$\<[\tau^1],[p]\>=(\tau^1\ci Tr)(p)$, where $p\in M_n(\asq)$, $p^2=p$, and
$Tr: M_n(\asq)\ra\asq$ is the usual matrix trace.
The following proposition establishes the pairing between the cyclic
cohomology class $[\tau^1]$ and the $K_0$-classes $[e_{-1}]$ and $[f_{-1}]$
of the left and right
projector matrix of bimodule $P_{-1}$ respectively:
\bpr\label{chern}
Let $\tau^1:\asq\ra k$ be the trace (\ref{tr}), and $e_{-1},f_{-1}$ 
the projectors given
in propositions \ref{ml} and~\ref{mr}. 
Then $(\tau^1\ci Tr)(e_{-1})=-1$ and $(\tau^1\ci Tr)(f_{-1})=1$.
\epr
\bpf
Taking advantage of (\ref{comm}) and (\ref{tr}), we get:
\[
(\tau^1\ci Tr)\pmatrix{\ha\hd & -\hb\ha \cr \hg\hd & -q\hb\hg \cr}
= \tau^1(1+(q^{-1}-q)\hb\hg)=\tau^1((q^{2}-1)\hz)=-1.
\]
Similarly,
\[
(\tau^1\ci Tr)\pmatrix{\hd\ha & \hd\hg \cr -\ha\hb & -q^{-1}\hb\hg \cr}
=1,
\]
as claimed.
\epf\\ \ \\
This computation is in agreement with the classical situation. Only there
the sign change of the Chern number when switching (by transpose) from
the left to right projector matrix is due to the anticommutativity of the
standard differential forms on manifolds. Here the sign change relies on
the noncommutativity of the algebra.

Since every free module can be represented in $K_0$ by the identity matrix,
we obtain that the pairing  of the cyclic cohomology class $[\tau^1]$ with
the $K_0$-class of any free \asq-module always vanishes:
\[
\<[\tau^1],[I]\>=\tau^1(\mbox{n})=0,~~ n\in\IN .
\]
Now, combining Proposition~\ref{chern} with Lemma~\ref{isol} yields:
\bco\label{noncl}
The \hge\ of the quantum principal Hopf fibration is {\em not} cleft.
\eco

\section*{Appendix}

In this appendix we provide a direct proof of non-cleftness of the quantum principal
Hopf fibration which is possible in the purely algebraic setting. This complements
our $K$-theoretic proof. Thus,
suppose that there exists a cleaving map 
$\Phi:k[z,z^{-1}]\ra\slq$. The existence of the convolution inverse $\Phi^{-1}$
entails  $\Phi(z)\Phi^{-1}(z)=\he(z)$, whence  $\Phi(z)$ must be invertible in~\slq.
The polynomial $\Phi(z)$ cannot be constant because then $\Phi(z)$ and $\Phi(1)=1$
would be linearly dependent, which contradicts the injectivity of $\Phi$
(see Section~1). Therefore to prove the non-cleftness it suffices  to show that all
invertible elements of \slq\ are non-zero numbers.

One can do it using the direct sum decomposition 
$\slq=\bigoplus_{m,n\in{\Bbb Z}}A[m,n]$, where
\[
A[m,n]=\{ p\in\slq\; |\;\pi(p\1)\ot p\2=z^m\ot p,\; 
p\1\ot\pi(p\2)=p\ot z^n\}
\]
(see~\cite[(1.10)]{mmnnu91}.) 
To be consistent with \cite{mmnnu91}, let us put now $k=\IC$. (See, however, bottom
of p.360 in~\cite{mmnnu91}.)
We know from \cite[p.363]{mmnnu91}
that we can write any element of \slq\ as a sum $\sum_{m,n}p_{m,n}(\zeta)e_{m,n}$
or $\sum_{k,l}e_{k,l}r_{k,l}(\zeta)$, where $\hz:=-q^{-1}\hb\hg$, 
$p_{m,n},r_{k,l}\in\IC[\zeta],\; 
e_{m,n}\in A[m,n]$. Assume now that 
$
\sum_{m,n}p_{m,n}(\zeta)e_{m,n} \sum_{k,l}e_{k,l}r_{k,l}(\zeta)=1.
$
Since both sums are finite, there exist indices
$m_+:=\mbox{max}\{m\in\IZ\;|\; p_{m,n}\neq 0\}$,
$n_+:=\mbox{max}\{n\in\IZ\;|\; p_{m_+,n}\neq 0\}$,
$m_-:=\mbox{min}\{m\in\IZ\;|\; p_{m,n}\neq 0\}$,
$n_-:=\mbox{min}\{n\in\IZ\;|\; p_{m_-,n}\neq 0\}$,
and similarly $k_+,k_-,l_+,l_-$.
We have
\beq\label{sum}
A[0,0]\ni e_{0,0}=1=\sum_{m,n}p_{m,n}(\zeta)e_{m,n} \sum_{k,l}e_{k,l}r_{k,l}(\zeta)
=\sum_{m,n,k,l}p_{m,n}(\zeta)s_{m,n,k,l}(\zeta)\tilde{r}_{k,l}(\zeta)e_{m+k,n+l}.
\eeq
Here $s_{m,n,k,l}(\zeta)e_{m+k,n+l}:=e_{m,n}e_{k,l}$ (see~\cite[p.363]{mmnnu91}),
and $\tilde{r}_{k,l}(\zeta)$ is obtained from $r_{k,l}(\zeta)$ by commuting it
over $e_{m+k,n+l}$, i.e., 
$e_{m+k,n+l}r_{k,l}(\zeta)=\tilde{r}_{k,l}(\zeta)e_{m+k,n+l}$. It follows from
the commutation relations (\ref{comm}) that the coefficients of $\tilde{r}_{k,l}$
are $q$ to some powers times the corresponding coefficients of~$r_{k,l}$. In 
particular, $r_{k,l}=0\Leftrightarrow\tilde{r}_{k,l}=0$. Since 
$p_{m_+,n_+}(\zeta)e_{m_+,n_+}$, $e_{k_+,l_+}r_{k_+,l_+}(\zeta)$ and
$p_{m_-,n_-}(\zeta)e_{m_-,n_-}$, $e_{k_-,l_-}r_{k_-,l_-}(\zeta)$ are the only 
terms that can contribute to the direct summand $A[m_++k_+,n_++l_+]$ and
$A[m_-+k_-,n_-+l_-]$ respectively, we can conclude from the equation~(\ref{sum})
that either $m_++k_+$, $n_++l_+$, $m_-+k_-$, $n_-+l_-$ are all zero, or else
$p_{m_{\pm},n_{\pm}}(\zeta)s_{m_{\pm},n_{\pm},k_{\pm},l_{\pm}}(\zeta)
\tilde{r}_{k_{\pm},l_{\pm}}(\zeta)e_{m_{\pm}+k_{\pm},n_{\pm}+l_{\pm}}=0$. From
\cite[p.363]{mmnnu91} we know, however, that $e_{m_{\pm}+k_{\pm},n_{\pm}+l_{\pm}}$
is a (left and right) basis of $A[m_{\pm}+k_{\pm},n_{\pm}+l_{\pm}]$ over 
$\IC[\hz]$. Also, using formulas $\ha^j\hd^j=\prod_{i=1}^j(1-q^{-2(i-1)}\hz)$, 
$\hd^j\ha^j=\prod_{i=1}^j(1-q^{2i}\hz)$ one can check that $e_{m,n}e_{k,l}\neq 0$, 
whence
$s_{m_{\pm},n_{\pm},k_{\pm},l_{\pm}}\neq 0$. Thus, as there are no zero divisors
in $\IC[\hz]$ and $r_{k,l}=0\Leftrightarrow\tilde{r}_{k,l}=0$, we can conclude
that $p_{m_{\pm},n_{\pm}}=0$ or $r_{k_{\pm},l_{\pm}}=0$. This, however, contradicts
the definition of $m_{\pm},n_{\pm},k_{\pm},l_{\pm}$. Therefore $m_{\pm}=-k_{\pm}$
and $n_{\pm}=-l_{\pm}$. Consequently, as $m_-\leq m_+$ and $k_-\leq k_+$, we have
$m_-=m_+=-k_+=-k_-$. Hence also $n_-=n_+=-l_+=-l_-$. Put $m_0=m_-=m_+$ and
$n_0=n_-=n_+$. It follows now that 
$\sum_{m,n}p_{m,n}(\zeta)e_{m,n}=p_{m_0,n_0}(\zeta)e_{m_0,n_0}$ and
$\sum_{k,l}e_{k,l}r_{k,l}(\zeta)=e_{-m_0,-n_0}r_{-m_0,-n_0}(\zeta)$. This way
(\ref{sum}) reduces to 
$p_{m_0,n_0}(\zeta)s_{m_0,n_0,-m_0,-n_0}(\hz)\tilde{r}_{-m_0,-n_0}(\zeta)=1$.
Hence all three of the above polynomials must be non-zero constants.
Using again \cite[p.363]{mmnnu91} and remembering that $\ha^j\hd^j$ and $\hd^j\ha^j$
are polynomials in \hz\ of degree $j$, we can infer that $m_0=0=n_0$. (Otherwise
$s_{m_0,n_0,-m_0,-n_0}$ is not of degree~$0$.) Consequently
$\sum_{m,n}p_{m,n}(\zeta)e_{m,n}=p_{0,0}(\zeta)$, 
$\sum_{k,l}e_{k,l}r_{k,l}(\zeta)=\tilde{r}_{0,0}(\zeta)=r_{0,0}(\zeta)$, and
$p_{0,0}$, $r_{0,0}$ are invertible constant polynomials, as needed.

\small

\section*{Acknowledgements} P.M.H. was partially supported 
by  the  NATO and CNR postdoctoral fellowships and KBN grants 2 P301 020 07
and 2 P03A 030 14. 
It is a pleasure to thank Max Karoubi and Giovanni Landi for very helpful 
discussions.


\begin{thebibliography}{MMNNU91}

\vspace*{-.mm}\bibitem[B-N72]{b-n72} Bourbaki, N.:~  
Commutative Algebra. Addison-Wesley 1972

\vspace*{-.mm}\bibitem[BM93]{bm93} Brzezi\'{n}ski, T.; Majid, S.:~
Quantum Group Gauge Theory on Quantum Spaces. 
Commun.\ Math.\ Phys.\ {\bf 157}  591--638 (1993); Erratum {\bf 167} 235 (1995)
\mbox{\sf(hep-th/9208007)}

\vspace*{-.mm}\bibitem[BM97]{bm97} Brzezi\'{n}ski, T.; Majid, S.:~
Quantum Differentials and the $q$-Monopole Revisited. 
To appear in Acta Applic.\ Math.\
\mbox{\sf(q-alg/9706021)}

\vspace*{-.mm}\bibitem[C-A94]{c-a94}  Connes, A.:~
Noncommutative Geometry. Academic Press 1994

\vspace*{-.mm}\bibitem[CQ95]{cq95} Cuntz, J.; D.Quillen, D.:~ 
 Algebra Extensions and Nonsingularity.
J.\ Amer.\ Math.\ Soc.\  {\bf 8} (2)  251--289 (1995)

\vspace*{-.mm}\bibitem[DH98]{dh98} D\c{a}browski, L.; Hajac, P.M.:~
Joint project. SISSA, Trieste, Italy, 1998 

\vspace*{-.mm}\bibitem[D-Y85]{d-y85} Doi, Y:~  
Algebras with total integrals.   
Commun.\ Alg.\ {\bf 13} 2137--2159 (1985)

\vspace*{-.mm}\bibitem[D-M96]{d-m96} Durdevic, M.:~  
Quantum Principal Bundles and Tannaka-Krein Duality Theory.    
Rep.\ Math.\ Phys.\ {\bf 38} (3) 313-324 (1996)
\mbox{\sf (q-alg/9507018)}


\vspace*{-.mm}\bibitem[H-PM96]{h-pm96} Hajac, P.M.:~  
Strong Connections on Quantum Principal Bundles.   
Commun.\ Math.\ Phys.\ {\bf 182} (3) 579--617 (1996)

\vspace*{-.mm}\bibitem[L-G97]{l-g97}  Landi, G.:~
An Introduction to Noncommutative Spaces and their Geometries. 
Springer-Verlag 1997

\vspace*{-.mm}\bibitem[L-JL97]{l-jl97}  Loday, J.-L.:~
Cyclic Homology. Springer 1997

\vspace*{-.mm}\bibitem[M-S95]{m-s95} Majid, S.:~
 Foundations of Quantum Group Theory.
Cambridge University Press 1995

\vspace*{-.mm}\bibitem[M-S97]{m-s97} Majid, S.:~
Some Remarks on Quantum and Braided Group Gauge Theory.
Banach Center Publications. {\bf 40} 336--349 (1997)
 
\vspace*{-.mm}\bibitem[MMNNU91]{mmnnu91} 
Masuda, T.; Mimachi, K.; Nakagami, Y.; Noumi, M.; Ueno, K.: 
Representations of the Quantum Group $SU\sb q(2)$ 
and the Little $q$-Jacobi Polynomials.
J.\ Funct.\ Anal.\ {\bf 99} 357--387 (1991)

\vspace*{-.mm}\bibitem[MNW91]{mnw91} 
Masuda, T.; , K.; Nakagami, Y.; Watanabe, J.:~ 
Noncommutative Differential Geometry on the Quantum Two Sphere of Podle\'s.
I: An Algebraic Viewpoint.
K-Theory {\bf 5} 151--175 (1991)

\vspace*{-.mm}\bibitem[P-P87]{p-p87} Podle\'{s}, P.:~  
Quantum Spheres.
Lett.\ Math.\ Phys.\ {\bf 14} 521--531 (1987)

\vspace*{-.mm}\bibitem[S-HJ90]{s-hj90} Schneider, H.-J.:~
Principal Homogenous Spaces for Arbitrary Hopf Algebras. 
Isr.\ J.\ Math.\ {\bf 72} (1--2) 167--195 (1990)

\vspace*{-.mm}\bibitem[S-HJ94]{s-hj94} Schneider, H.J.:~ 
Hopf Galois Extensions, Crossed Products, and Clifford Theory. 
In: Bergen, J., Montgomery, S. (eds.)~  
Advances in Hopf Algebras. Lecture Notes in Pure and Applied 
Mathematics. Marcel Dekker, Inc. {\bf 158} 267--297 (1994)

\end{thebibliography}
\end{document}